\documentclass[11pt]{article}
\usepackage[]{amsmath,amssymb}

\newtheorem{theorem}{Theorem}[section]
\newtheorem{lemma}[theorem]{Lemma}
\newtheorem{proposition}[theorem]{Proposition}
\newtheorem{corollary}[theorem]{Corollary}
\newtheorem{exAux}[theorem]{Example}

\newtheorem{Note}[theorem]{Note}
\newenvironment{note}{\begin{Note} \rm}{\end{Note}}
\newtheorem{Def}[theorem]{Definition}
\newenvironment{definition}{\begin{Def} \rm}{\end{Def}}
\newtheorem{Rem}[theorem]{Remark}
\newenvironment{remark}{\begin{Rem} \rm}{\end{Rem}}
\newtheorem{Ass}[theorem]{Assumption}

\newenvironment{proof}{\medskip\noindent{\bf Proof.\ }}{\qed\medskip}
\newenvironment{proofof}[1]{\medskip\noindent{\bf Proof  of {#1}.\ 
}}{\qed\medskip}
\newcommand{\qed}{\hfill\mbox{$\Box$\qquad\qquad}}
\newcommand{\Mat}[1]{\text{\rm Mat}_{#1}(\mathbb{K})}

\begin{document}
\thispagestyle{empty}
\vspace{1cm}

\noindent
\begin{center}
\LARGE \bf
The determinant of $AA^*-A^*A$ for \\ a Leonard pair $A,A^*$
\end{center}

\medskip\noindent
\begin{center}
\Large
Kazumasa Nomura and Paul Terwilliger
\end{center}

\smallskip

\begin{quote}
\small 
\begin{center} {\bf Abstract}
\end{center}
Let $\mathbb{K}$ denote a field, and let $V$ denote a vector
space over $\mathbb{K}$ with finite positive dimension.
We consider a pair of linear transformations
$A:V \to V$ and $A^* : V \to V$
that satisfy (i), (ii) below:
\begin{itemize}
\item[(i)] There exists a basis for $V$ with respect to which
the matrix representing $A$ is irreducible tridiagonal and the
matrix representing $A^*$ is diagonal.
\item[(ii)] There exists a basis for $V$ with respect to which
the matrix representing $A^*$ is irreducible tridiagonal and the
matrix representing $A$ is diagonal.
\end{itemize}
We call such a pair a {\em Leonard pair} on $V$.
In this paper we investigate the commutator $AA^*-A^*A$.
Our results are as follows.
Abbreviate $d=\dim V-1$ and 
first assume $d$ is odd. We show $AA^*-A^*A$ is invertible and
display several attractive formulae for the determinant.
Next assume $d$ is even.
We show that the null space of $AA^*-A^*A$ has dimension $1$.
We display a nonzero vector in this null space.
We express this vector as a sum of eigenvectors for $A$ and
as a sum of eigenvectors for $A^*$.
\end{quote}

\section{Introduction}

Throughout this paper, $\mathbb{K}$ will denote a field
and $V$ will denote a vector space over $\mathbb{K}$ with finite
positive dimension.

We begin by recalling the notion of a Leonard pair.
We will use the following notation.
A square matrix $X$ is called {\em tridiagonal}
whenever each nonzero entry lies on either the diagonal, the subdiagonal,
or the superdiagonal. Assume $X$ is tridiagonal.
Then $X$ is called {\em irreducible}
whenever each entry on the subdiagonal is nonzero and each entry on
the superdiagonal is nonzero.

\medskip

\begin{definition}   \cite{T:Leonard}      \label{def:LP}
By a {\em Leonard pair} on $V$, we mean an ordered pair of 
linear transformations $A:V \to V$ and $A^* : V \to V$
that satisfy the following two conditions:
\begin{itemize}
\item[(i)] There exists a basis for $V$ with respect to which
the matrix representing $A$ is irreducible tridiagonal and the
matrix representing $A^*$ is diagonal.
\item[(ii)] There exists a basis for $V$ with respect to which
the matrix representing $A^*$ is irreducible tridiagonal and the
matrix representing $A$ is diagonal.
\end{itemize}
\end{definition}

\medskip

\begin{note}
It is a common notational convention to use $A^*$ to represent the
conjugate-transpose of $A$. We are {\em not} using this convention.
In a Leonard pair $A$, $A^*$ the linear transformations $A$ and
$A^*$ are arbitrary subject to (i) and (ii) above.
\end{note}

\medskip
We refer the reader to 
\cite{H},
\cite{N:aw}, 
\cite{NT:balanced}, \cite{NT:formula},
\cite{P}, \cite{T:sub1}, \cite{T:sub3}, \cite{T:Leonard},
\cite{T:24points}, \cite{T:conform}, \cite{T:intro},
\cite{T:intro2}, \cite{T:split}, \cite{T:array}, \cite{T:qRacah},
\cite{T:survey}, \cite{TV}, \cite{V}
for background on Leonard pairs.
We especially recommend the survey \cite{T:survey}.
See \cite{AC}, \cite{AC2}, \cite{ITT}, \cite{IT:shape},
\cite{IT:uqsl2hat}, \cite{IT:non-nilpotent}, \cite{ITW:equitable}, 
\cite{N:refine}, \cite{N:height1},
\cite{T:qSerre}, \cite{T:Kac-Moody} for related topics.

\medskip
For the rest of this paper let $A$, $A^*$ denote a Leonard
pair on $V$. For notational convenience we define $d = \dim V -1$.
We are going to investigate the commutator $AA^* - A^*A$.
It turns out the behavior of this commutator depends on the parity of
$d$. First assume $d$ is odd. We show $AA^*-A^*A$ is invertible and
display several attractive formulae for the determinant.
Next assume $d$ is even.
We show that the null space of $AA^*-A^*A$ has dimension $1$.
We display a nonzero vector in this null space.
We express this vector as a sum of eigenvectors for $A$ and
as a sum of eigenvectors for $A^*$.

The rest of this section is devoted to giving formal statements
of our results.
We will use the following notation.
Let $\Mat{d+1}$ denote the $\mathbb{K}$-algebra consisting of all $d+1$ by
$d+1$ matrices that have entries in $\mathbb{K}$. 
We index the rows and  columns by $0, 1, \ldots, d$. 
Let $u_0,u_1,\ldots,u_d$ denote a basis for $V$. 
For a linear transformation $X:V \to V$ and for 
$Y \in \Mat{d+1}$, we say 
{\em $Y$ represents $X$ with respect to} $u_0, u_1, \ldots, u_d$
whenever $X u_j = \sum_{i=0}^d Y_{ij}u_i$ for $0 \leq j \leq d$.
Also, by the {\em null space} of $X$ we mean $\{ v\in V\,|\, Xv=0\}$.
We recall that $X$ is invertible if and only if the null space 
of $X$ is zero.
We now state our first main result.
\medskip

\begin{theorem}   \label{thm:rank}
The following hold.
\begin{itemize}
\item[(i)] Suppose $d$ is odd. Then $AA^*-A^*A$ is invertible.
\item[(ii)] Suppose $d$ is even. Then the null space
of $AA^*-A^*A$ has dimension $1$.
\end{itemize}
\end{theorem}

\bigskip
Before we state our next two main theorems we make some comments.
We fix a basis $v^*_0, v^*_1, \ldots, v^*_d$ for $V$ that
satisfies Definition \ref{def:LP}(i). 
Observe that with respect to this basis the matrices that represent $A$, $A^*$ 
take the form
\begin{equation}     \label{eq:AAs}
A : 
\begin{pmatrix}
 a_0 & b_0 & & & &  \text{\bf 0} \\
 c_1 & a_1 & b_1 \\
 & c_2 & \cdot & \cdot \\
 & & \cdot & \cdot & \cdot \\
 & & & \cdot & \cdot & b_{d-1} \\
 \text{\bf 0} & & &  & c_d & a_d
\end{pmatrix},
\qquad
A^* : 
\begin{pmatrix}
 \theta^*_0 & & & & & \text{\bf 0} \\
 & \theta^*_1 \\
 & & \theta^*_2 \\
 & & & \cdot \\
 & & & & \cdot \\
 \text{\bf 0} &  & & & & \theta^*_d
\end{pmatrix}
\end{equation}
where $b_{i-1}c_{i} \neq 0$ for $1 \leq i \leq d$.
We also fix a basis $v_0, v_1, \ldots, v_d$ for $V$ that
satisfies Definition \ref{def:LP}(ii). 
With respect to this basis the matrices that represent $A$, $A^*$ 
take the form
\begin{equation}     \label{eq:AAs2}
A : 
\begin{pmatrix}
 \theta_0 & & & & & \text{\bf 0} \\
 & \theta_1 \\
 & & \theta_2 \\
 & & & \cdot \\
 & & & & \cdot \\
 \text{\bf 0} &  & & & & \theta_d
\end{pmatrix},
\qquad
A^* :
\begin{pmatrix}
 a^*_0 & b^*_0 & & & &  \text{\bf 0} \\
 c^*_1 & a^*_1 & b^*_1 \\
 & c_2 & \cdot & \cdot \\
 & & \cdot & \cdot & \cdot \\
 & & & \cdot & \cdot & b^*_{d-1} \\
 \text{\bf 0} & & &  & c^*_d & a^*_d
\end{pmatrix}
\end{equation}
where $b^*_{i-1}c^*_{i} \neq 0$ for $1 \leq i \leq d$.
Observe that $\theta_0, \theta_1, \ldots, \theta_d$
(respectively $\theta^*_0, \theta^*_1, \ldots, \theta^*_d$)
are the eigenvalues of $A$ (respectively $A^*$).
It is known \cite[Lemma 1.3]{T:Leonard} that 
$\theta_i \neq \theta_j$, $\theta^*_i \neq \theta^*_j$
if $i\neq j$ for $0 \leq i,j \leq d$.
We now state our second and third main results.

\medskip

\begin{theorem}           \label{thm:det1}      \samepage
Suppose $d$ is odd. Then
\begin{eqnarray} 
\text{\rm det} (AA^*-A^*A)
 &=& \prod_{\stackrel{1 \leq i \leq d}{\text{ $i$ odd}}}
      b_{i-1}c_{i}(\theta^*_{i-1}-\theta^*_i)^2,
        \label{eq:det1} \\
\text{\rm det} (AA^*-A^*A)
 &=& \prod_{\stackrel{1 \leq i \leq d}{\text{ $i$ odd}}}
      b^*_{i-1}c^*_{i}(\theta_{i-1}-\theta_i)^2.
        \label{eq:det1s}
\end{eqnarray}
\end{theorem}

\medskip

\begin{theorem}          \label{thm:span}  
Suppose $d$ is even.
\begin{itemize}
\item[(i)] 
The null space of $AA^*-A^*A$ is spanned by
$\sum_{k=0}^{d} \gamma_k v^*_k$, 
where $\gamma_k=0$ if $k$ is odd, and
\begin{equation}           \label{eq:lambda}
   \gamma_k = 
        \prod_{ \stackrel{1 \leq i \leq k-1}{ \text{$i$ odd} } }
        \frac{c_i(\theta^*_{i-1}-\theta^*_{i})}
             {b_i (\theta^*_{i}-\theta^*_{i+1})}
\end{equation}
if $k$ is even.
\item[(ii)]
The null space of $AA^*-A^*A$ is spanned by
$\sum_{k=0}^{d} \gamma^*_k v_k$, 
where $\gamma^*_k=0$ if $k$ is odd, and
\begin{equation}           \label{eq:lambdas}
   \gamma^*_k = 
        \prod_{ \stackrel{1 \leq i \leq k-1}{ \text{$i$ odd} } }
        \frac{c^*_i(\theta_{i-1}-\theta_{i})}
             {b^*_i (\theta_{i}-\theta_{i+1})}
\end{equation}
if $k$ is even.
\end{itemize}
\end{theorem}

\medskip

\begin{remark}
Theorems \ref{thm:rank} and \ref{thm:span} give an answer to
a problem by the second author \cite[Section 36]{T:survey}.
\end{remark}

\medskip

In order to state our next result we recall a few facts.

\medskip

\begin{lemma}     \cite[Theorem 1.9]{T:Leonard}   \label{lem:beta}
The expressions
\begin{equation}           \label{eq:beta}
    \frac{\theta_{i-2}-\theta_{i+1}}
         {\theta_{i-1}-\theta_{i}},
   \qquad \qquad
    \frac{\theta^*_{i-2}-\theta^*_{i+1}}
         {\theta^*_{i-1}-\theta^*_{i}}
\end{equation}
are equal and independent of $i$ for $2 \leq i \leq d-1$.
\end{lemma}

\medskip
Using Lemma \ref{lem:beta} we define a scalar $q$ as follows.

\medskip

\begin{definition}          \label{def:beta}
For $d \geq 3$ let $\beta$ denote the scalar in $\mathbb{K}$ such that 
$\beta+1$ is the common value of (\ref{eq:beta}).
For $d \leq 2$ let $\beta$ denote any scalar in $\mathbb{K}$.
Let $\overline{\mathbb{K}}$ denote the algebraic closure of $\mathbb{K}$.
Let $q$ denote a nonzero scalar in  $\overline{\mathbb{K}}$ such that
$\beta=q^2+q^{-2}$. 
\end{definition}

\medskip

We recall some notation.

\medskip

\begin{definition}         \label{def:qi}
For an integer $n > 0$ we define
\begin{equation}         \label{eq:qint}
    [n]_q = q^{n-1}+q^{n-3} + \cdots + q^{1-n}.
\end{equation}
We observe:
\begin{itemize}
\item[(i)] If $q^2 \neq 1$ then
\begin{equation}     \label{eq:qint1}
     [n]_q = \frac{q^n-q^{-n}}{q-q^{-1}}.
\end{equation}
\item[(ii)] If $q=1$ then
\begin{equation}     \label{eq:qint2}
  [n]_q=n.
\end{equation}
\item[(iii)] If $q=-1$ then
\begin{equation}     \label{eq:qint3}
   [n]_q = (-1)^{n+1} n.
\end{equation}
\end{itemize}
\end{definition}

\medskip

We mention here a technical result for later use.
We will show
\begin{equation}        \label{eq:qintnonzero}
   [i]_q \neq 0 \;\;\; \text{ if $i$ is odd}
   \qquad (1 \leq i \leq d).
\end{equation}

\medskip
We recall some parameters.
By \cite[Theorem 3.2]{T:Leonard}
there exists a sequence of nonzero scalars 
$\varphi_1, \varphi_2, \ldots, \varphi_d$ in $\mathbb{K}$ and 
there exists a basis for $V$ with respect to which the matrices representing 
$A$, $A^*$ are
\[
A :
\begin{pmatrix}
 \theta_0 & & & & \text{\bf 0} \\
 1 & \theta_1 \\
 & 1 & \theta_2 \\
 & &  \cdot & \cdot \\
 & &  & \cdot & \cdot \\
 \text{\bf 0} & & & & 1 & \theta_d
\end{pmatrix},
\qquad
A^* :
\begin{pmatrix}
 \theta^*_0 & \varphi_1 & & & & \text{\bf 0} \\
 & \theta^*_1 & \varphi_2 \\
 & & \theta^*_2 & \cdot \\
 & & & \cdot & \cdot \\
 & & & & \cdot & \varphi_d \\
 \text{\bf 0} & & & &  & \theta^*_d
\end{pmatrix}.
\]
The sequence $\varphi_1, \varphi_2, \ldots, \varphi_d$
is uniquely determined by the ordering
$(\theta_0, \theta_1, \ldots, \theta_d$;
$\theta^*_0, \theta^*_1, \ldots,$ $\theta^*_d)$.
We call the sequence $\varphi_1, \varphi_2, \ldots, \varphi_d$
the {\em first split sequence} with respect to the ordering
$(\theta_0, \theta_1, \ldots, \theta_d$; 
$\theta^*_0, \theta^*_1, \ldots, \theta^*_d)$.
Let $\phi_1, \phi_2, \ldots, \phi_d$ denote the first split
sequence with respect to the ordering
$(\theta_d, \theta_{d-1}, \ldots, \theta_0$; 
$\theta^*_0, \theta^*_1, \ldots, \theta^*_d)$.
We call the sequence $\phi_1, \phi_2, \ldots, \phi_d$ the
{\em second split sequence} with respect to the ordering
$(\theta_0, \theta_1, \ldots, \theta_d$; 
$\theta^*_0, \theta^*_1, \ldots, \theta^*_d)$.
We now state our final main result.

\medskip

\begin{theorem}        \label{thm:det2}      \samepage
Suppose $d$ is odd. Then 
\begin{equation}       \label{eq:det2}
   \text{\rm det} (AA^*-A^*A)
     = (-1)^{(d+1)/2} \prod_{\stackrel{1 \leq i \leq d}{\text{ $i$ odd}}}
          \frac{\, \varphi_i \phi_i \,}{[i]_q^{\,2}}.
\end{equation}
\end{theorem}

\medskip
\begin{remark}
The denominator of (\ref{eq:det2}) is nonzero by (\ref{eq:qintnonzero}).
\end{remark}

\begin{remark}
Theorem \ref{thm:det2} was conjectured
by the second author \cite[Section 36]{T:survey}.
\end{remark}

\section{Proof of Theorems \ref{thm:rank}, \ref{thm:det1} and  \ref{thm:span}}

\begin{lemma}     \label{lem:B}
Let $B \in \Mat{d+1}$ denote the matrix that represents $AA^*-A^*A$
with respect to the basis $v^*_0, v^*_1, \ldots, v^*_d$.
Then
\begin{itemize}
\item[(i)] The $(i,i-1)$-entry of $B$ is 
$c_i(\theta^*_{i-1}-\theta^*_{i})$ for $1 \leq i \leq d$.
\item[(ii)] The $(i-1,i)$-entry of $B$ is 
$b_{i-1}(\theta^*_{i}-\theta^*_{i-1})$ for $1 \leq i \leq d$.
\item[(iii)] All other entries of $B$ are $0$.
\end{itemize}
\end{lemma}

\begin{proof}
Obtained by routine computation using (\ref{eq:AAs}).
\end{proof}

\medskip

\begin{proofof}{Theorem \ref{thm:det1}}
We first prove (\ref{eq:det1}).
Let the matrix $B$ be as in Lemma \ref{lem:B}. 
Observe that $B$ is tridiagonal with all diagonal entries zero.
For $0 \leq r \leq d$  let $B_r$ denote the submatrix of $B$
obtained by taking rows $0,1,\ldots,r$ 
and columns $0,1,\ldots,r$.
Then the determinants of $B_1,B_3,\ldots,B_d$ satisfy the following 
well-known recursion \cite[p.~28]{HJ}:
\[
 \det(B_1)=b_0 c_1(\theta^*_0 - \theta^*_1)^2,
\]
\[
   \det(B_r) = b_{r-1}c_{r}(\theta^*_{r-1}-\theta^*_{r})^2 \det(B_{r-2})
   \qquad(3 \leq r \leq d, \text{ $r$ odd}).
\]
Solving this recursion we find
\begin{equation}      \label{eq:det1aux}
 \det(B_r) = 
     \prod_{\stackrel{1 \leq i \leq r}{i \text{ odd}}} 
          b_{i-1}c_{i} (\theta^*_{i-1}-\theta^*_{i})^2 
   \qquad (1 \leq r \leq d, \;\;\text {$r$  odd}).
\end{equation}
Setting $r=d$ in (\ref{eq:det1aux}) we obtain (\ref{eq:det1}). 
The proof of (\ref{eq:det1s}) is similar.
\end{proofof}

\medskip

\begin{proofof}{Theorem \ref{thm:span}(i)}
Let the matrix $B$ be as in Lemma \ref{lem:B}.
Define a vector $v=(\gamma_0,\gamma_1,\ldots,\gamma_d)^t$
where $t$ denotes the transpose, and $\gamma_0, \gamma_1, \ldots, \gamma_d$ 
are from the statement of the theorem.
It suffices to show that $v$ spans the null space of $B$.
By matrix multiplication we find $Bv=0$,
so $v$ is contained in the null space of $B$.
Let $w$ denote any vector in the null space of $B$. 
We show $w$ is a scalar multiple of $v$.
For notational convenience write $w=(w_0,w_1,\ldots,w_d)^t$.
Multiplying out $Bw=0$ we routinely obtain
the recursion
\[
    b_0(\theta^*_1 - \theta^*_0) w_1 = 0,
\]
\[
    c_{i}(\theta^*_{i-1}-\theta^*_{i})w_{i-1}
  + b_{i}(\theta^*_{i+1}-\theta^*_{i})w_{i+1} = 0
    \qquad (1 \leq i \leq d-1),
\]
\[
    c_{d}(\theta^*_{d-1}-\theta^*_{d})w_{d-1} = 0.
\]
Solving this recursion we find $w_k=w_0 \gamma_k$ for $0 \leq k \leq d$.
Therefore $w=w_0 v$. 
We have now shown that $w$ is a scalar multiple of $v$ and the result follows.
\end{proofof}

\medskip

\begin{proofof}{Theorem \ref{thm:span}(ii)}
Similar to the proof of Theorem \ref{thm:span}(i).
\end{proofof}

\medskip

\begin{proofof}{Theorem \ref{thm:rank}}
Immediate from Theorems \ref{thm:det1} and \ref{thm:span}.
\end{proofof}

\section{The proof of Theorem \ref{thm:det2}, part I}

We now turn to the proof of Theorem \ref{thm:det2}.
We will use the following notation.
Let $\lambda$ denote an indeterminate and let 
$\mathbb{K}[\lambda]$ denote the $\mathbb{K}$-algebra consisting of all polynomials 
in $\lambda$ that have coefficients in $\mathbb{K}$.

\medskip

\begin{definition}
For $0 \leq i \leq d$ let $\tau^*_i$,
$\eta^*_i$ denote the following polynomials in $\mathbb{K}[\lambda]$:
\begin{eqnarray}
 \tau^*_i &=& (\lambda-\theta^*_0)(\lambda-\theta^*_1) \cdots(\lambda-\theta^*_{i-1}),
    \label{eq:deftaus} \\
 \eta^*_i &=& 
   (\lambda-\theta^*_d)(\lambda-\theta^*_{d-1})\cdots (\lambda-\theta^*_{d-i+1}).
      \label{eq:defetas}
\end{eqnarray}
We observe that each of $\tau^*_i$, $\eta^*_i$ is monic with degree $i$.
\end{definition}

\medskip

\begin{theorem}   \cite[Lemma 7.2, Theorem 23.7]{T:survey} \label{thm:bc}
For $1 \leq i \leq d$ we have
\begin{equation}        \label{eq:xi}
    b_{i-1}c_{i} = \varphi_i \phi_i \:
       \frac{\tau^*_{i-1}(\theta^*_{i-1})\eta^*_{d-i}(\theta^*_i)}
            {\tau^*_{i}(\theta^*_i)\eta^*_{d-i+1}(\theta^*_{i-1})}.
\end{equation}
\end{theorem}

\medskip
We now assume that $d$ is odd and evaluate
(\ref{eq:det1}) using (\ref{eq:xi}).  We find that
$\det(AA^*-A^*A)$ is equal to
\[
   \prod_{\stackrel{1 \leq i \leq d}{i \text{ odd}}} \varphi_i \phi_i
\]
times
\begin{equation}          \label{eq:left}
\prod_{\stackrel{1 \leq i \leq d}{i \text{ odd}}}
    (\theta^*_{i-1}-\theta^*_{i})^2
     \frac{\tau^*_{i-1}(\theta^*_{i-1})\eta^*_{d-i}(\theta^*_{i})}
          {\tau^*_{i}(\theta^*_{i})\eta^*_{d-i+1}(\theta^*_{i-1})}.
\end{equation}
We now evaluate (\ref{eq:left}).

\medskip

\begin{lemma}         \label{lem:1}
Suppose $d$ is odd. Then (\ref{eq:left}) is equal to
\begin{equation}        \label{eq:1}
  (-1)^{m+1} \Psi^2,
\end{equation}
where $m=(d-1)/2$ and
\[
  \Psi = 
\prod_{0 \leq \ell<k \leq m}
     \frac{\theta^*_{2\ell+1}-\theta^*_{2k}}{\theta^*_{2\ell}-\theta^*_{2k+1}}.
\]
\end{lemma}

\begin{proof}
For an integer $i$ define $s(i)=(-1)^i$.
Using (\ref{eq:defetas}) we find
\[
\prod_{\stackrel{1 \leq i \leq d}{i \text{ odd}}}
     \frac{\eta^*_{d-i}(\theta^*_{i})}
          {\eta^*_{d-i+1}(\theta^*_{i-1})}
= \prod_{0 \leq i < j \leq d} (\theta^*_{i}-\theta^*_{j})^{s(i+1)}.
\]
Similarly using (\ref{eq:deftaus}) we find
\[
\prod_{\stackrel{1 \leq i \leq d}{i \text{ odd}}}
     \frac{\tau^*_{i-1}(\theta^*_{i-1})}
          {\tau^*_{i}(\theta^*_{i})}
= (-1)^{m+1} \prod_{0 \leq i < j \leq d} (\theta^*_{i}-\theta^*_{j})^{s(j)}.
\]
Evaluating (\ref{eq:left}) using these equations we routinely obtain the result.
\end{proof}

\section{Some comments}

In order to prove Theorem \ref{thm:det2} we will evaluate (\ref{eq:1})
further using Lemma \ref{lem:beta}.
There are some technical aspects involved which we will deal with
in this section.

\medskip

\begin{lemma}   \cite[Lemma 9.3]{T:Leonard}    \label{lem:char}
Assume $d \geq 3$. Then with reference to Definition \ref{def:beta} 
the following (i)--(iv) hold.
\begin{itemize}
\item[(i)] Suppose $\beta \neq 2$, $\beta \neq -2$. 
Then $q^{2i} \neq 1$ for $1 \leq i \leq d$.
\item[(ii)] Suppose $\beta=2$ and $\text{\rm Char}(\mathbb{K})=p > 2$.
Then $d<p$.
\item[(iii)] Suppose $\beta=-2$ and $\text{\rm Char}(\mathbb{K})=p > 2$.
Then $d < 2p$.
\item[(iv)] Suppose $\beta=0$ and $\text{\rm Char}(\mathbb{K})=2$. Then
$d = 3$.
\end{itemize}
\end{lemma}

\medskip

\begin{lemma}             \label{lem:added}
Referring to Definition \ref{def:qi}, assume that $n$ is odd and $q^2=-1$. 
Then $[n]_q=(-1)^{(n-1)/2}$.
\end{lemma}

\begin{proof}
Routine using line (\ref{eq:qint1}).
\end{proof}

\medskip

\begin{corollary}    \label{cor:nonzero}
With reference to Definitions \ref{def:beta} and \ref{def:qi}, 
we have $[i]_q \neq 0$ for $i$ odd, $(1 \leq i \leq d)$.
\end{corollary}

\begin{proof}
Assume $d\geq 3$; otherwise the result holds since $[1]_q=1$.
Let the integer $i$ be given and assume $i$ is odd.
We consider three cases. 
First assume $\beta \neq 2$, $\beta \neq -2$.
Then the result holds by Lemma \ref{lem:char}(i) and (\ref{eq:qint1}).
Next assume $\beta=2$ and $\text{\rm Char}(\mathbb{K})\neq 2$.
Using $\beta=2$ and $q^2+q^{-2}=\beta$
we find $q^2=1$. Now $[i]_q = i$ by (\ref{eq:qint2}) or (\ref{eq:qint3}) and 
since $i$ is odd.
Each of $1,2,\ldots, d$ is nonzero in $\mathbb{K}$ by
Lemma \ref{lem:char}(ii) so $[i]_q \neq 0$.
Next assume $\beta =-2$.
Using $q^2+q^{-2}=\beta$ we find $q^2=-1$, 
so $[i]_q = (-1)^{(i-1)/2}$ by Lemma \ref{lem:added}.
In particular  $[i]_q \neq 0$ as desired.
\end{proof}

\medskip

\begin{lemma}    \cite[Lemma 9.4]{T:Leonard}    \label{lem:frac}
Assume $d \geq 3$.
Pick any integers $i$, $j$, $r$, $s$ $(0 \leq i,j,r,s \leq d)$ and
assume $i+j=r+s$, $i \neq j$, $r \neq s$. 
Then with reference to Definition \ref{def:beta} the following (i)--(iv) hold.
\begin{itemize}
\item[(i)] Suppose $\beta \neq 2$, $\beta\neq -2$. Then
\[
     \frac{\theta^*_i -\theta^*_j}{\theta^*_r - \theta^*_s}
     = \frac{q^{2i} - q^{2j}}{q^{2r}-q^{2s}}.
\]
\item[(ii)] Suppose $\beta=2$ and $\text{\rm Char}(\mathbb{K}) \neq 2$. Then
\[
     \frac{\theta^*_i -\theta^*_j}{\theta^*_r - \theta^*_s}
     = \frac{i-j}{r-s}.
\]
\item[(iii)]
Suppose $\beta=-2$ and $\text{\rm Char}(\mathbb{K}) \neq 2$.
If $r+s$ is even, then
\[
     \frac{\theta^*_i -\theta^*_j}{\theta^*_r - \theta^*_s}
     = (-1)^{i+r} \frac{i-j}{r-s}.
\]
If $r+s$ is odd, then
\[
     \frac{\theta^*_i -\theta^*_j}{\theta^*_r - \theta^*_s}
     = (-1)^{i+r}.
\]
\item[(iv)]
Suppose $\beta=0$ and $\text{\rm Char}(\mathbb{K}) = 2$. Then
\[
     \frac{\theta^*_i -\theta^*_j}{\theta^*_r - \theta^*_s}
     = 1.
\]
\end{itemize}
In the above formulae all denominators are nonzero by Lemma \ref{lem:char}.
\end{lemma}

\medskip

\begin{corollary}      \label{cor:1}
Assume $d \geq 3$.
Pick any integers $i$, $j$ $(0 \leq i<j \leq d)$ and
$r$, $s$ $(0 \leq r<s \leq d)$.
Assume $i+j=r+s$ and this common value is odd.
Then with reference to Definitions \ref{def:beta} and \ref{def:qi},
\begin{equation}         \label{eq:frac}
\frac{\theta^*_i - \theta^*_j}
     {\theta^*_r - \theta^*_s}
= \frac{[j-i]_q}{[s-r]_q}.
\end{equation}
\end{corollary}

\begin{proof}
In each case of Lemma \ref{lem:frac} we routinely express the
result using Definition \ref{def:qi} and Lemma \ref{lem:added}.
\end{proof}

\section{Proof of Theorem \ref{thm:det2}, part II}

In this section we complete the proof of Theorem \ref{thm:det2}.
Our argument is based on the following proposition.

\medskip

\begin{proposition}     \label{prop:2}
Assume $d$ is odd. Then the expression $\Psi$ from
Lemma \ref{lem:1} satisfies
\begin{equation}    \label{eq:2}
\Psi = \prod_{\stackrel{1 \leq i \leq d}{i \text{ odd}}} \frac{1}{[i]_q}.
\end{equation}
\end{proposition}

\begin{proof}
We may assume $d \geq 3$; otherwise the result holds since $[1]_q=1$.
Now we have
\begin{eqnarray*}
    \Psi &=& \prod_{0 \leq \ell<k\leq m}  
                \frac{\theta^*_{2\ell+1}-\theta^*_{2k}}
                     {\theta^*_{2\ell}-\theta^*_{2k+1}}
\\
        &=& \prod_{k=0}^m  \prod_{\ell=0}^{k-1} 
                \frac{\theta^*_{2\ell+1}-\theta^*_{2k}}
                     {\theta^*_{2\ell}-\theta^*_{2k+1}}
\\
        &=& \prod_{k=0}^m  \prod_{\ell=0}^{k-1} 
                \frac{[2k-2\ell-1]_q}{[2k-2\ell+1]_q} 
		                     \qquad \text{(by Corollary \ref{cor:1})}
\\
        &=& \prod_{k=0}^m  \frac{1}{[2k+1]_q}
\\
        &=& \prod_{\stackrel{1 \leq i \leq d}{i \text{ odd}}}
               \frac{1}{[i]_q}.
\end{eqnarray*}
\end{proof}

\medskip

\begin{proofof}{Theorem \ref{thm:det2}}
Immediate from Lemma \ref{lem:1}, Proposition \ref{prop:2}, and 
the comment after Theorem \ref{thm:bc}.
\end{proofof}

%
%
%
%
%
%
%
%

\bigskip

\bibliographystyle{plain}

\begin{thebibliography}{1}

\bibitem{AC}
H.~Alnajjar  and B.~Curtin.
\newblock
A family of tridiagonal pairs.
\newblock {\em
Linear Algebra Appl.}
{\bf 390}
(2004)
369--384.

\bibitem{AC2}
H.~Alnajjar  and B.~Curtin.
\newblock
A family of tridiagonal pairs related to
the quantum affine algebra 
$U\sb q(\widehat{\mathfrak{sl}}\sb 2)$.
\newblock {\em
Electron. J. Linear Algebra}
{\bf 13} (2005) 1--9. 


\bibitem{H}
B.~ Hartwig. 
Three mutually adjacent Leonard pairs. 
{\em Linear Algebra Appl.} 
To appear.

\bibitem{HJ}
R.~A.~Horn and C.~R.~Johnson.
Matrix Analysis.
Cambridge University Press, Cambridge, 1990.


\bibitem{ITT}
T.~Ito, K.~Tanabe, and P.~Terwilliger.
\newblock Some algebra related to ${P}$- and ${Q}$-polynomial association
  schemes,  in:
\newblock {\em Codes and Association Schemes (Piscataway NJ, 1999)}, Amer.
Math. Soc., Providence RI, 2001, pp.
     167--192; 
{\tt arXiv:math.CO/0406556}.

\bibitem{IT:shape}
T.~Ito and P.~Terwilliger.
\newblock The shape of a tridiagonal pair.
\newblock {\em J. Pure Appl. Algebra} {\bf 188} (2004) 145--160;
{\tt arXiv:math.QA/0304244}.

\bibitem{IT:uqsl2hat}
T.~Ito and P.~Terwilliger.
\newblock {Tridiagonal pairs and the quantum affine 
algebra
$U_q({\widehat{sl}}_2)$.}
\newblock {\em Ramanujan J.}, accepted; 
{\tt arXiv:math.QA/0310042}.

\bibitem{IT:non-nilpotent}
T.~Ito and P.~Terwilliger.
\newblock
Two non-nilpotent linear transformations that satisfy the cubic $q$-Serre relations.
{\em J. Algebra}. Submitted;
{\tt arXiv:math.QA/0508398}.

\bibitem{ITW:equitable}
T.~Ito, P.~Terwilliger and C.~Weng.
\newblock {The quantum algebra $U_q(sl_2)$ and its equitable presentation.}
{\em J. Algebra}. To appear;
{\tt arXiv:math.math.QA/0507477}.

\bibitem{N:aw}
K.~Nomura.
\newblock
Tridiagonal pairs and the {A}skey-{W}ilson relations.
\newblock {\em Linear Algebra Appl.}
{\bf 397} (2005) 99--106.

\bibitem{N:refine}
K.~Nomura.
\newblock A refinement of the split decomposition of
a tridiagonal pair.
\newblock {\em Linear Algebra Appl.}
{\bf 403} (2005) 1--23.

\bibitem{N:height1}
K.~Nomura.
\newblock
Tridiagonal pairs of height one.
\newblock {\em Linear Algebra Appl.}
{\bf 403} (2005) 118--142.

\bibitem{NT:balanced}
K.~Nomura and P.~Terwilliger.
\newblock
Balanced Leonard pairs.
\newblock {\em Linear Algebra Appl.}
Submitted;
{\tt arXiv:math.RA/0506219}. 

\bibitem{NT:formula}
K.~Nomura and P.~Terwilliger.
\newblock
Some formulae involving the split sequences of a Leonard pair.
{\em Linear Algebra Appl.}
To appear;
{\tt arXiv:math.RA/0508407}

\bibitem{P}
A.~A.~Pascasio.
\newblock     
On the multiplicities of the primitive idempotents of a
 {$Q$}-polynomial distance-regular graph.
\newblock{ \em
European J. Combin.}
{\bf 23}
(2002)
1073--1078.

\bibitem{T:sub1}
P.~Terwilliger.
\newblock The subconstituent algebra of an association scheme I. 
\newblock {\em J. Algebraic Combin.} {\bf 1} (1992) 363--388.
   
\bibitem{T:sub3}
P.~Terwilliger.
\newblock The subconstituent algebra of an association scheme III.
\newblock{\em
J. Algebraic Combin. }
{\bf 2}  (1993) 177--210.

\bibitem{T:Leonard}
P.~Terwilliger.
\newblock Two linear transformations each tridiagonal with respect to an
  eigenbasis of the other.
  \newblock {\em Linear Algebra Appl.}  {\bf 330} (2001) 149--203;
{\tt arXiv:math.RA/0406555}.

\bibitem{T:qSerre}
  P.~Terwilliger.
  \newblock Two relations that generalize the $q$-Serre relations and the
  Dolan-Grady relations. In
  \newblock {\em  Physics and
  Combinatorics 1999 (Nagoya)}, 377--398, World Scientific Publishing,
   River Edge, NJ, 2001; 
{\tt arXiv:math.QA/0307016}.

\bibitem{T:24points}
   P.~Terwilliger.
   \newblock  Leonard pairs from 24 points of view.
   \newblock {\em Rocky Mountain J. Math.} {\bf 32}(2) (2002) 827--888;
{\tt arXiv:math.RA/0406577}.

\bibitem{T:conform}
   P.~Terwilliger.
   \newblock Two linear transformations each tridiagonal with respect to an
     eigenbasis of the other; the $TD$-$D$ and the $LB$-$UB$ canonical form.
\newblock {\em J. Algebra}, to appear;
{\tt arXiv:math.RA/0304077}.

\bibitem{T:intro}
    P.\ Terwilliger.
    \newblock Introduction to Leonard pairs.
    \newblock {OPSFA Rome 2001}.
    \newblock{\em J. Comput. Appl. Math.} {\bf 153}(2) (2003)
    463--475.

\bibitem{T:intro2}
P.\ Terwilliger.
Introduction to {L}eonard pairs and
  {L}eonard systems.
 {\em S\=uri\-kaiseki\-kenky\=usho} {\em K\=oky\=uroku},
 (1109) 67--79, 1999.   Algebraic combinatorics  (Kyoto, 1999).

\bibitem{T:split}
P.~Terwilliger.
\newblock Two linear transformations each tridiagonal with respect to an
  eigenbasis of the other; comments on the split decomposition.
\newblock {\em 
 J. Comput. Appl. Math.} {\bf 178} (2005) 437--452;
{\tt arXiv:math.RA/0306290}.

 \bibitem{T:array}
 P.~Terwilliger.
 \newblock Two linear transformations each tridiagonal with respect to an
   eigenbasis of the other; comments on the parameter array.
\newblock {\em
Des. Codes Cryptogr.}  {\bf 34}  (2005) 307--332;
{\tt arXiv:math.RA/0306291}.

\bibitem{T:qRacah}
P.~Terwilliger.
\newblock Leonard pairs and the $q$-Racah polynomials.
\newblock {\em Linear Algebra Appl.} {\bf 387} (2004) 235--276;
{\tt arXiv:math.QA/0306301}.

\bibitem{T:survey}
P.~Terwilliger.
\newblock
Two linear transformations each tridiagonal with respect to an eigenbasis
of the other; an algebraic approach to the Askey scheme of
orthogonal polynomials.
\newblock Lecture notes for the summer school on orthogonal polynomials
and special functions. Universidad Carlos III de Madrid, Leganes, Spain.
July 8--July 18, 2004;
{\tt arXiv:math.QA/0408390}. 

\bibitem{T:Kac-Moody}
P.~Terwilliger.
\newblock
The equitable presentation for the quantum group $U_q(g)$ 
associated with a symmetrizable Kac-Moody algebra $g$.
{\em J. Algebra}. Submitted; 
{\tt arXiv:math.QA/0507478}.


\bibitem{TV}
P.~Terwilliger and R.~Vidunas.
\newblock Leonard pairs and the Askey-Wilson relations.
\newblock {\em J. Algebra Appl.} {\bf 3} (2004) 411--426;
{\tt arXiv:math.QA/0305356}.

\bibitem{V}
R.~Vidunas.
\newblock
Normalized Leonard pairs and Askey-Wilson relations. Preprint.
\hfil\break
{\tt 
arXiv:math.RA/0505041}.

\end{thebibliography}

\bigskip\bigskip\noindent
Kazumasa Nomura\\
College of Liberal Arts and Sciences\\
Tokyo Medical and Dental University\\
Kohnodai\\
Ichikawa, 272-0827 Japan\\
email: nomura.las@tmd.ac.jp

\bigskip\noindent
Paul Terwilliger\\
Department of Mathematics\\
University of Wisconsin\\
Van Vleck Hall\\
480 Lincoln drive\\
Madison, WI 53706-1388 USA\\
email: terwilli@math.wisc.edu

\bigskip\noindent
{\bf Keywords.}
Leonard pair, Terwilliger algebra, Askey scheme,
$q$-Racah polynomial.

\noindent
{\bf 2000 Mathematics Subject Classification}.
Primary: 15A15.
Secondary: 05E35, 05E30, 33C45, 33D45.

\end{document}